\documentclass[10pt,openbib]{article}

\usepackage{amsfonts,amsthm}
\usepackage{amsmath}
\usepackage{amssymb}
\usepackage{graphicx}

\newcommand{\R}{\mathbb{R}}
\newcommand{\C}{\mathbb{C}}

\newcommand{\onehalf}{\mbox{${\scriptstyle \frac{1}{2}\, }$}}
\newcommand{\dee}{\mathop{\! \, \rm d \!}\nolimits}
\newcommand{\comp}{\, \raisebox{2pt}{$\scriptstyle\circ \, $}}

\newcommand{\setrule}{\, \rule[-4pt]{.5pt}{13pt}\, }

\newcommand{\spann}{\mathop{\rm span}\nolimits}
\newcommand{\lefthook}{\mbox{$\, \rule{8pt}{.5pt}\rule{.5pt}{6pt}\, \, $}}

\newcommand{\ttfrac}[2]{\mbox{$\frac{{\scriptstyle #1}}{{\scriptstyle #2}}$}}

\newcommand{\circdot}{\raisebox{2pt}{\tiny$\, \bigodot$}}

\allowdisplaybreaks

\begin{document}

\begin{center}
{\Large \textbf{On affine Riemann surfaces}} \\ 
\vspace{.05in}
Richard Cushman
\end{center}
\addtocounter{footnote}{1}
\footnotetext{printed: \today}

\begin{abstract}
We show that the universal covering space of a connected component of a regular level set of a smooth complex valued function on ${\C}^2$, which is a smooth affine Riemann surface, is ${\R }^2$. This implies that the orbit 
space of the action of the covering group on ${\R }^2$ is the original affine Riemann surface. 
\end{abstract} 

\section{Basic result}
%%%%%%%%%%%%%%

Let 
\begin{displaymath}
F : {\C}^2 \rightarrow \C : (z,w) \mapsto u + \mathrm{i}\, v = \mathrm{Re}\, F + \mathrm{i}\, \mathrm{Im}\, F
\end{displaymath}
be a smooth function. Let 
$X_F$ be the holomorphic Hamiltonian vector field on $({\C}^2 , \dee z \wedge \dee w )$ corresponding to $F$, 
that is, $X_F \lefthook (\dee z \wedge \dee w) = \dee F$. 
On ${\C }^2 = {\R }^4 = (\mathrm{Re}\, z, \mathrm{Im}\, z, \mathrm{Re}\, w, \mathrm{Im}\, w)$ with real symplectic form 
$\Omega = \mathrm{Re}(\dee z \wedge \dee w)$ we have real Hamiltonian vector fields $X_u$ and $X_v$. Then 
\begin{align*}
(X_u + \mathrm{i} \, X_v) \lefthook \mathrm{Re}(\dee z \wedge \dee w) & =  \dee u + \mathrm{i}\, \dee v = \dee F 
= X_F \lefthook (\dee z \wedge \dee w) \\
&\hspace{-1.5in} = ( \mathrm{Re}\, X_F + \mathrm{i}\, \mathrm{Im}\, X_F) \lefthook 
\big( \mathrm{Re}\, (\dee z \wedge \dee w) 
+ \mathrm{i}\, \mathrm{Im}\, (\dee z \wedge \dee w) \big) \\
&\hspace{-1.5in}  = \big( \mathrm{Re}\, X_F + \mathrm{i}\, \mathrm{Im}\, X_F \big) \lefthook \mathrm{Re}(\dee z \wedge \dee w) + \big( - \mathrm{Im}\, X_F +\mathrm{i}\, \mathrm{Re}\, X_F \big) \lefthook \mathrm{Im}(\dee z \wedge \dee w) 
\end{align*}
So 
\begin{displaymath}
(X_u + \mathrm{i} \, X_v) \lefthook \mathrm{Re}(\dee z \wedge \dee w) = 
\big( \mathrm{Re}\, X_F + \mathrm{i}\, \mathrm{Im}\, X_F \big) \lefthook \mathrm{Re}(\dee z \wedge \dee w), 
\end{displaymath}
since the $2$-forms $\mathrm{Re}(\dee z \wedge \dee w)$ and $\mathrm{Im}(\dee z \wedge \dee w)$ are 
linearly independent. This implies 
\begin{displaymath}
X_{\mathrm{Re}\, F} = X_u = \mathrm{Re}\, X_F \quad \mathrm{and} \quad X_{\mathrm{Im}\, F} = 
X_v = \mathrm{Im}\, X_F,
\end{displaymath} 
since $\mathrm{Re}(\dee z \wedge \dee w)$ is nondegenerate. \medskip 

\noindent \textbf{Proposition 1.1} Let $S$ be a connected component of  $F^{-1}(c)$, where $c \in \C $ is a regular value of $F$, which lies in its image. Then the universal covering space of $S$ is ${\R }^2$. \medskip 

\noindent \textbf{Proof.} $S$ is a smooth $1$ dimensional complex manifold, which we 
assume is connected. Our argument constructs coordinates on the universal covering space of $S$. 
We begin. For every $(z,w) \in S$ the complex tangent space to $S$ at $(z,w)$ is 
$\ker \dee F (z,w)$, where 
\begin{displaymath}
(0,0) \ne \dee F (z,w) = \dee u_{|S}(z,w) + \mathrm{i}\, \dee v_{|S}(z,w) = (\dee u + \mathrm{i} \, \dee v)_{|S}(z,w) 
\end{displaymath}%
for every $(z,w) \in S$. Thus the nonzero vector field $X_F = (X_u + \mathrm{i}\, X_v)_{|S}$ spans the 
complex tangent space of $S$ at each point of $S$.  Because 
$X_F$ is nonzero on $S$, the real vector fields ${X_u}_{|S}$ and 
${X_{v}}_{|S}$ are linearly independent at each point of $S$. To see this we argue as follows. 
Suppose that the real vector fields $X_{u_{|S}}$ and $X_{v_{|S}}$ are linearly dependent at some 
point $(z,w) \in S$. Then ${\spann }_{\R } \{ X_{u_{|S}}(z,w), X_{v_{|S}}(z,w) \}$ has real dimension $1$. 
Thus $(X_{u_{|S}} + \mathrm{i}\, X_{v_{|S}})(z,w)$ does not span the complex tangent space to $S$ at 
$(z,w)$, which is a contradiction. \medskip 

\vspace{-.15in}Consider the $2$-form ${\Omega }_{|S}$ on $S$. Since $\Omega $ is closed, it follows that 
${\Omega }_{|S}$ is closed. Because ${X_u}_{|S}$ and ${X_v}_{|S}$ are linearly independent 
vector fields on $S$ and $\Omega $ is nondegenerate on ${\R}^4$, it follows that ${\Omega }_{|S}$ 
is nondegenerate on ${\spann }_{\R } \{ {X_u}_{|S}(z,w),$ \linebreak 
${X_v}_{|S}(z,w) \} $ for every $(z,w) \in S$. To see this from $(X_u + \mathrm{i}\, X_v) \lefthook 
\Omega = \dee F$ and the fact that $\dee F \ne (0,0)$ on $S$ we get $\Omega (X_u , X_v) \ne 0$ on $S$. 
Hence ${\Omega }_{|S}$ is a symplectic form on $S$. \medskip 

Let $M$ be the universal covering space of $S$ with covering mapping $\rho : M \rightarrow S$. Because 
$\rho $ is a local diffeomorphism, the $2$-form $\omega = {\rho }^{\ast} ({\Omega }_{|S})$ on $M$ is 
symplectic. Consider  the smooth functions $U = {\rho }^{\ast }(u_{|S})$ and $V = {\rho }^{\ast }(v_{|S})$ on 
$(M, \omega )$. The corresponding Hamiltonian vector fields $X_U$ and $X_V$ on $(M, \omega )$ are given by 
$\dee U = X_U \lefthook \omega $ and $\dee V = X_V \lefthook \omega $. Since 
\begin{align*}
X_U \lefthook \omega & = \dee U = \dee \, ({\rho }^{\ast}u_{|S}) = {\rho }^{\ast }(\dee u_{|S}) 
 = {\rho }^{\ast }(X_{u_{|S}} \lefthook {\Omega }_{|S} ) \\
 &  = {\rho }^{\ast}(X_{u_{|S}}) \lefthook {\rho }^{\ast }({\Omega }_{|S}) = {\rho }^{\ast }(X_{u_{|S}}) \lefthook \omega , 
\end{align*}
it follows that $X_U = {\rho }^{\ast }(X_{u_{|S}})$, because $\omega $ is nondegenerate. Similarly, 
$X_V = {\rho }^{\ast }(X_{v_{|S}})$. Since $\rho $ is a local diffeomorphism and the vector fields 
$X_{u_{|S}}$ and $X_{v_{|S}}$ are linearly independent at each point of $S$, the vector fields 
$X_U$ and $X_V$ are linearly independent at each point of $M$. Thus the $1$-forms 
$\dee U$ and $\dee V$ on $M$ are linearly independent at each point of $M$, because $\omega $ is 
nondegenerate. So the vector fields $\frac{\partial }{\partial U}$ and $\frac{\partial }{\partial V}$ 
are linearly independent at each point of $M$. \medskip 

Consider the nonzero $2$-form $\varpi = \dee V \wedge \dee U$ on $M$. 
Since $M$ is $2$-dimensional, the de Rham cohomology group of $2$-forms on $M$ has dimension $1$. 
Thus $\varpi = a \omega $ for some nonzero real number $a$.\footnote{We compute $a$ as follows. 
Let $D \subseteq {\R }^2 $ be the unit disk in $({\R }^2, \varpi =\dee V \wedge \dee U)$ with Euclidean 
inner product. Orient $D$ so that its boundary 
is traversed clockwise. Then $\pi =\int_D\varpi = a \int_D \omega $, that is, $a = \pi /\int_D \omega $.} Because 
$\{ \frac{\partial }{\partial U}, \frac{\partial }{\partial V} \}$ is a basis of the tangent space of $M$ at 
each point of $M$, we may write $X_U = A \frac{\partial }{\partial U} + B\frac{\partial }{\partial V}$. Then 
\begin{displaymath}
\dee U = X_U \lefthook \omega = \ttfrac{1}{a} X_U \lefthook \varpi = \ttfrac{1}{a}(B \dee U - A \dee V), 
\end{displaymath}
which implies $X_U = a \frac{\partial }{\partial V}$. A similiar argument shows that 
$X_V = -a \frac{\partial }{\partial U}$. \medskip 

The pair of functions $(U,V)$ are coordinates on $M$, since the vector fields $X_U = 
a \frac{\partial }{\partial V}$ and $X_V = -a \frac{\partial }{\partial U}$ are linearly independent at each point of $M$ and commute. This latter assertion follows because 
\begin{align*}
\{ u, v \} & = L_{X_v}u = L_{X_{\mathrm{Im}\, F}}(\mathrm{Re}\, F) 
= L_{\ttfrac{1}{2\mathrm{i}}(X_{F - \mathrm{i}\, F})}\onehalf (F +\mathrm{i}\, F) \\
& = \ttfrac{1}{4i}[ L_{X_F}F + \mathrm{i} \, L_{X_F}F - \mathrm{i}\, L_{X_F}F + L_{X_F}F ]  = 0 , 
\end{align*}
implies $[X_v, X_u] = X_{\{ u,v \} } = 0$. From 
\begin{displaymath}
T\rho \, [X_U,X_V] = [X_u|S, X_v|S ] \comp \rho  = {[X_u, X_v]}_{|S} \comp \rho = 0, 
\end{displaymath}
we get $[X_U,X_V] =0$, because $\rho $ is a local diffeomorphism. Thus 
we may identify $M$ with ${\R }^2$. \hfill $\square $ \medskip 

\noindent \textbf{Corollary 1.1A}(Bates and Cushman \cite{bates-cushman}). The image of the linear 
flow of the vector field $X_{U+ \mathrm{i}\, V}$ on $\C$ under the covering map $\rho $ is 
the flow of the vector field $X_F$ on $S$.  \medskip

\noindent \textbf{Proof.} The flow of $X_{U + \mathrm{i}\, V}$ on $\C$ is $U(t) + \mathrm{i}\, V(t) = 
\big( U(0) + \mathrm{i}at \big) +\big( \mathrm{i}V(0) -at \big)$, since $X_U = a \frac{\partial }{\partial V}$ and 
$X_V = -a \frac{\partial }{\partial U}$. Hence an integral curve of $X_{U + \mathrm{i}\, V}$ starting at 
$U(0) + \mathrm{i}\, V(0)$ is $t \mapsto \big( U(0) + \mathrm{i}\, V(0) \big) + a( -t + \mathrm{i}t)$, which is 
a straight line in $\C$. Thus the flow of $X_{U + \mathrm{i} V}$ is linear. Since 
\begin{align*}
T\rho X_{U + \mathrm{i}\, V} & = T\rho (X_U + \mathrm{i}\, X_V) = T\rho X_U + \mathrm{i}\, T\rho X_V \\
& = X_{u_{|S}} \comp \rho + \mathrm{i}\, X_{v_{|S}} \comp \rho = X_{(u+\mathrm{i}v)_{|S}} \comp \rho =
{X_{F}}_{|S} \comp \rho , 
\end{align*}
the image of the flow of $X_{U+ \mathrm{i}\, V}$ under the covering map $\rho $ 
is the flow of $X_F$. \hfill $\square $ \medskip 

Define a Riemannian metric $E$ on ${\R }^2$ by $\mathrm{E} = \ttfrac{1}{a^2} \dee U \circdot \dee U 
+ \ttfrac{1}{a^2} \dee V \circdot \dee V$. Since $\mathrm{E }(X_U, X_U) = 1 = \mathrm{E} (X_V,X_V)$ and 
$\mathrm{E} (X_U, X_V) =0$, we find that $\mathrm{E}$ is the Euclidean inner product on 
$T_{(U,V)}{\R }^2 = {\R }^2$ for every $(U,V) \in {\R }^2$. The metric $\mathrm{E}$ is flat, since 
it is indendent of $(U,V) \in {\R }^2$. Let $G$ be the group of covering transformations of $S$. 
Then $G$ is a discrete subgroup of the two dimensional Euclidean group. $G$ acts properly 
on ${\R }^2$. Since each element of $G$ leaves no point of ${\R }^2$ fixed, we obtain the \medskip 

\noindent \textbf{Corollary 1.1B.} The orbit space ${\R }^2/G$ of the action on $S$ of the covering group $G$ 
on the universal covering space ${\R }^2$ of the affine Riemann surface $S$ is diffeomorphic to $S$.

\section{Example$^{\, \scriptscriptstyle 3}$}
\addtocounter{footnote}{1}
\footnotetext{See Cushman \cite{cushman}.}
%%%%%%%%%%%

Let 
\begin{equation}
F: {\C}^2 \rightarrow \C: (z,w) \mapsto w^2 + z^6. 
\label{eq-ex-one}
\end{equation}
Then $1$ is a regular value of $F$, since $(0,0) = \dee F(z,w) =(6z^5, 2w)$ if and only if 
$z = w =0$. But $(0,0) \notin F^{-1}(1) = S$. Thus $S$ is a smooth affine Riemann surface. 
Let $\pi : {\C}^2 \rightarrow \C : (z,w) \mapsto z$. Then ${\pi }_{|S}: S \subseteq {\C}^2 \rightarrow \C$ 
is a branched covering map of $S$ with branch points $B = \{ (z_k = {\mathrm{e}}^{2\pi \mathrm{i}k/6}, 0) \in S \setrule \,  
\mbox{for $k=0,1, \ldots , 5$} \}$ and branch values $V = \{ z_k \setrule \, k=0,1, \ldots , 5 \}$. 
The map ${\pi }_{|S}$ is smooth on $S \setminus B$ with image $\C \setminus V$. The sheets $S_{\ell }$ of the 
branched covering map ${\pi }_{|S}$ are defined by $w_{\ell } = {\mathrm{e}}^{2\pi \mathrm{i} \ell /2}(1-z^6)^{1/2}$ for 
$\ell = 0,1$, where $z \in \C $, that is, $S_{\ell }$ is a connected component of 
$({\pi }_{|S})^{-1}(\C ) = \coprod_{\ell = 0,1}S_{\ell }$. \medskip 

Let $\rho : {\R }^2 \rightarrow S$ be the universal covering map of $S$. The sheets of the covering map 
$\rho $ are ${\Sigma }_{\ell } = {\rho }^{-1}(S_{\ell })$ for $\ell =0,1$. The group $G$ of covering transformations 
of $S$ is the collection of isometries of $({\R }^2, \mathrm{E})$, where $\mathrm{E}$ is the Euclidean inner 
product on ${\R }^2$, which permute the sheets ${\Sigma }_{\ell }$ of $\rho $. Consider the group $G'$ of 
diffeomorphisms of $S$ generated by the transformations 
\begin{align*}
&\mathcal{R}: S\subseteq {\C}^2 \rightarrow S\subseteq {\C}^2: (z,w) \mapsto ({\mathrm{e}}^{2\pi \mathrm{i}/6}z, w) 
\intertext{and} 
&\mathcal{U}:S\subseteq {\C}^2 \rightarrow S\subseteq {\C}^2: (z,w) \mapsto (\overline{z}, \overline{w}) . 
\end{align*}
Since ${\mathcal{R}}^6 = {\mathcal{U}}^2 = \mathrm{id}$ and $\mathcal{R}\mathcal{U} = \mathcal{U}{\mathcal{R}}^{-1}$, 
the group $G'$ is isomorphic to the dihedral group on $6$ letters.\footnote{The group $G'$ is also 
generated by the reflections $\{ R^kU, \, k=0,1, \ldots 5 \setrule \, R^6 = U^2= \mathrm{id}  \} $. Thus $G'$ is the Weyl group of the complex simple Lie algebra ${\mathbf{A}}_5$.} 
Because $\mathcal{R}(S_{\ell }) = S_{\ell }$ for 
$\ell =0,1$ and $\mathcal{U}(S_0) = S_1$, the map $\mathcal{R}$ induces the identity permutation of the sheets of 
the covering map $\rho $; while the map $\mathcal{U}$ transposes the sheets of $\rho $. Thus 
$\mathcal{R}$ and $\mathcal{U}$ generate the covering group $G$. \medskip 

We want to describe the action of $G$, as a subgroup of the Euclidean group of $({\R }^2, \mathrm{E})$. \medskip 

We will need some preliminary results. Let 
\begin{equation}
f: \C \setminus V \rightarrow \C : z \mapsto \int^z_0 \frac{1}{2w} \, \dee z , 
\label{eq-ex-two}
\end{equation}
where $w = \sqrt{1 - z^6}$. Then $f$ is a local diffeomorphism, because $\dee f = \frac{1}{2w} \dee z$ is 
nonvanishing on $\C \setminus V$. We have \medskip 

\noindent \textbf{Proposition 2.1} Up to a coordinate transformation $\lambda : \C \rightarrow \C$, the map 
\begin{equation}
\delta : S \subseteq {\C}^2 \rightarrow \C : (z,w) \mapsto \zeta = \alpha (f \comp {\pi }_{|S}) (z,w) , 
\label{eq-ex-three}
\end{equation}
where $\alpha = \sqrt{2} {\mathrm{e}}^{3\pi \mathrm{i}/4}$, is a right inverse of the 
universal covering map $\rho $, that is, $\rho \comp \lambda \comp \delta = {\mathrm{id}}_S$. \medskip 

To prove proposition 2.1 we need \medskip 

\noindent \textbf{Lemma 2.2} The image under the map $\delta $ (\ref{eq-ex-three}) of an 
integral curve of the vector field $(X_F)_{|S}$ on $S$ is an integral curve of the vector field 
$\alpha \frac{\partial }{\partial \zeta }$ on $\C$. \medskip 

\noindent \textbf{Proof.} It suffices to show that for every $(z,w) \in S$ 
\begin{equation} 
T_{(z,w)} \delta \, X_F(z,w) = 
\alpha \frac{\partial }{\partial \zeta } \rule[-10pt]{.5pt}{22pt}\raisebox{-9pt}{$\, \scriptstyle{\zeta = \delta (z,w)}$}
\hspace{-25pt}. 
\label{eq-ex-four}
\end{equation}
This we do as follows. Using the definition of the map ${\pi }_{|S}$ and the vector field $(X_F)_{|S} = 
2w \frac{\partial }{\partial z} - 6w^5\frac{\partial }{\partial w}$, for every $(z,w) \in S$ we get  
\begin{displaymath}
T_{(z,w)}{\pi }_{|S} \, X_F(z,w) = T_{(z,w)}{\pi }_{|S}(2w \frac{\partial }{\partial z} - 6w^5\frac{\partial }{\partial w}) =
2w \frac{\partial }{\partial z} . 
\end{displaymath}
By definition of the function $f$ (\ref{eq-ex-two}) we have $\dee f = \frac{1}{2w} \dee z$, which implies 
$T_zf \big( 2w \frac{\partial }{\partial z} \big) $ $= \frac{\partial }{\partial \zeta }$. Thus for every 
$(z,w) \in S$ 
\begin{align*}
T_{(z,w)}\delta \, X_F(z,w) &= \alpha T_zf \Big( T_{(z,w)}{\pi }_{|S} \big( X_F(z,w) \big) \Big) = 
\alpha \frac{\partial }{\partial \zeta }, 
\end{align*} 
which establishes equation (\ref{eq-ex-four}). \hfill $\square $ \medskip

\noindent \textbf{Corollary 2.2A} The map $\delta $ (\ref{eq-ex-three}) is a local diffeomorphism. \medskip 

\noindent \textbf{Proof.} This follows from equation (\ref{eq-ex-four}), which shows that the 
tangent map of $\delta $ is injective at each point of $S$. \hfill $\square $ \medskip 

\noindent \textbf{Proof of proposition 2.1} Let 
$U + \mathrm{i}V = {\rho }^{\ast }(\mathrm{Re}\, F) + \mathrm{i}\, {\rho }^{\ast }(\mathrm{Im}\, F) $. 
By proposition 1.1, $U+\mathrm{i}\, V$ is a coordinate on $\C$. Define the diffeomorphism  
\begin{displaymath}
\lambda : \C \rightarrow \C: \zeta \mapsto U + \mathrm{i}\, V
\end{displaymath}
by requiring ${\lambda }_{\ast }\big( \alpha \frac{\partial }{\partial \zeta } \big) = X_U + \mathrm{i}\, X_V$, that is, 
set $U = \lambda (\mathrm{Re}\,\zeta )$ and $V = \lambda (\mathrm{Im}\, \zeta )$. By 
construction we have $\alpha \frac{\partial }{\partial \zeta } = {\lambda}^{\ast }{\rho }^{\ast } \big( (X_F)_{|S} \big)$, 
see the proof of proposition 1.1. By equation (\ref{eq-ex-four}) we have $\alpha \frac{\partial }{\partial \zeta } = 
{\delta }_{\ast }\big( (X_F)_{|S} \big)$. Thus ${\delta }_{\ast } = {\lambda}^{\ast }{\rho }^{\ast }$, which implies 
$\rho \comp \lambda \comp \delta = {\mathrm{id}}_S$. To see this suppose that $\rho \comp \lambda \comp 
\delta \ne {\mathrm{id}}_S$. Then ${\delta }^{\ast } \comp (\rho \comp \lambda )^{\ast } \ne \mathrm{id}_{TS}$. 
Hence ${\lambda }^{\ast }{\rho }^{\ast } \ne {\delta }_{\ast }$, which is a contradiction. \hfill $\square $ \medskip 

Let
\begin{equation}
R: \C \rightarrow \C : z \mapsto {\mathrm{e}}^{2\pi \mathrm{i}/6}z.
\label{eq-ex-five}
\end{equation}
Then $f(Rz) = Rf(z)$,  where $f$ is the function defined in (\ref{eq-ex-two}). To see this we compute. 
\begin{align*}
f(Rz) & = \int^{Rz}_0 \frac{\dee \xi}{2w(\xi)}, \, \, \mbox{where $w(\xi ) = \sqrt{1-{\xi }^6}$} \\
& = \int^z_0 \frac{R \dee z}{2 w(z)}, \, \, \mbox{using $\xi = Rz$ and $w(Rz) = w(z)$} \\
& = Rf(z). 
\end{align*}
Thus the image under $f$ (\ref{eq-ex-two}) of the closed equilateral triangle 
\begin{displaymath}
T' = \{ z = r' {\mathrm{e}}^{\mathrm{i}  {\theta}' } \in \C \setrule \, 0 \le r' \le 1 \, \, 
\& \, \, 0 \le {\theta}' \le 2\pi /6 \} 
\end{displaymath}
with vertex at the origin and one edge of length $1$ along the real axis is the equilateral 
triangle 
\begin{displaymath}
T = f(T') = \{ \zeta = r {\mathrm{e}}^{\mathrm{i} \theta } \in \C \setrule \, 0 \le r \le C \, \, 
\& \, \, 2\pi /6 \le \theta \le 4\pi /6 \} = RT', 
\end{displaymath}
where $C=\int^1_0 \frac{\dee z}{\sqrt {1-z^6}}$. Hence $f$ maps a regular hexagon into another. In particular, 
it sends the closed regular hexagon $H'$ with center at the origin $O$ and 
edge length $1$ onto the regular hexagon $H$ with center at $O$ and edge length $C$. Since 
$H'$ is simply connected and is contained in the unit disk $\{ | z | \le 1 \}$, the complex square root 
$\sqrt{1-z^6}$ is single valued for all $z \in H'$. Thus $H'$ is the image under 
${\pi }_{|S}$ of a domain $\mathcal{D} \subset S$, which is contained in some sheet $S_{{\ell }'}$ 
of the covering map $\rho $ of $S$. \medskip 

Let 
\begin{equation}
U: \C \rightarrow \C : z \mapsto \overline{z}. 
\label{eq-ex-six}
\end{equation}
The regular hexagon $H$ is invariant under the action of the group $\mathcal{G}$, generated 
be the rotation $R$ and the reflection $RU$ in the diagonal of $H$, which is 
an edge of the triangle $T$ that is not the real axis. The map $\delta $ (\ref{eq-ex-three}) intertwines the 
action of the group $G'$ generated by $\mathcal{R}$ and 
$\mathcal{RU}$ on $S$ with the action of the $\mathcal{G}$ on $H$. Thus the domain $\mathcal{D}$ 
contains a fundamental domain of the action of the covering group $G$ on ${\R }^2$. \medskip 

Let $\mathcal{T}$ be the abelian group generated by the translations 
\begin{displaymath}
{\tau }_k : \C \rightarrow \C : z \mapsto z + u_k,
\, \, \mbox{for $k=0,1, \ldots , 5$.}
\end{displaymath}
Here $u_k = \sqrt{3}C\, {\mathrm{e}}^{2\pi \mathrm{i}(1/12+ k/6)}$, which is perpendicular to an edge of 
the equilateral triangle $R^k(T)$ that lies on the boundary of the hexagon $H$. The action of 
$\mathcal{T}$ on $\C$ has fundamental domain $H$. To see this recall that in \cite{cushman} it is shown that 
\begin{displaymath}
\bigcup_{n \ge 0}\hspace{-6pt} \bigcup_{\hspace{10pt}{\ell }_1 + \cdots + {\ell}_k =n} \hspace{-15pt} {\tau }^{{\ell }_1}_1 \raisebox{-1pt}{$\comp$} \cdots \raisebox{-1pt}{$\comp$} {\tau }^{{\ell }_k}_k(K) = \C, 
\end{displaymath}
where $K$ is the closed stellated hexagon formed by placing an equilateral triangle of edge length $C$ on 
each bounding edge of $H$. But 
\begin{displaymath}
K = H \cup \bigcup^5_{k=0} {\tau }_k(R^{(4+k) \! \bmod 6}T). 
\end{displaymath}
So $H$ is the fundamental domain of the $\mathcal{T}$ action on $\C$. Because 
applying an element of $G'$ to the domain $\mathcal{D} \subseteq S$ gives a 
domain whose boundary has a nonempty intersection with the boundary of 
$\mathcal{D}$, it follows that under the mapping $\delta $ (\ref{eq-ex-three}) the 
corresponding element of the group of motions in $\C$ sends the hexagon $H$ 
to a hexagon which has an edge in common with $H$. Thus this group of motions 
is the group $\mathcal{T}$. Because the mapping $\delta $ intertwines the 
$G'$ action on $S$ with the $\mathcal{T}$ action on $\C$ and sends the domain $\mathcal{D} \subseteq S_{{\ell }'}$ 
diffeomorphically onto $H$, it follows that $\mathcal{D}$ is a fundamental domain for the action 
of $G'$ on $S$. Consider $\lambda (H)$, which is a regular hexagon with center 
at the origin, since the coordinate change $\lambda $ maps straight lines to straight lines. From proposition 2.1 we deduce that $\lambda (H)$ is a fundamental domain for the action of the 
covering group $G$ on $\C = {\R }^2$ of the affine Riemann surface $S$. Hence 
$S = {\R }^2/\mathcal{T}$. \hfill $\square $

\end{document}